\documentclass{article}

\usepackage{amssymb}
\usepackage{amsthm}

\newtheorem{thm}{Theorem}[section]
\newtheorem{prop}[thm]{Proposition}

\newtheorem{cor}[thm]{Corollary}

\newtheorem{defn}{Definition}

\begin{document}

\title{Structures in Familiar Classes Which Have Scott Rank $\omega_1^{CK}$}

\author{W.\ Calvert, S.\ S.\ Goncharov, and J.\ F.\ Knight}

\maketitle

\begin{abstract}

There are familiar examples of computable structures having various
computable Scott ranks.  There are also familiar structures, such as
the Harrison ordering, which have Scott rank $\omega_1^{CK}+1$.  Makkai
\cite{M} produced a structure of Scott rank $\omega_1^{CK}$, which can
be made computable \cite{K-Y}, and simplified so that it is just a
tree \cite{C-K-Y}.  In the present paper, we show that there are
further computable structures of Scott rank $\omega_1^{CK}$ in the
following classes:  undirected graphs, fields of any characteristic,
and linear orderings.  The new examples share with the Harrison
ordering, and the tree in \cite{C-K-Y}, a strong approximability
property. 

\end{abstract}   

\section{Introduction}

In this section, we recall some definitions and earlier results.  Scott rank is a measure of model-theoretic complexity.  The notion comes from the Scott Isomorphism Theorem (see \cite{S}, or \cite{Ke}).

\begin{thm} [Scott Isomorphism Theorem] 
\label{thm1.1} 

For each countable structure $\mathcal{A}$ (for a countable language $L$) there is an $L_{\omega_1\omega}$ sentence whose countable models 
are just the isomorphic copies of $\mathcal{A}$.

\end{thm}  

In the proof, Scott assigned countable ordinals to tuples in $\mathcal{A}$, and to $\mathcal{A}$ itself.
There are several different definitions of \emph{Scott rank} in use.    
We begin with a family of equivalence relations.  

\begin{defn} 

Let $\overline{a}$, $\overline{b}$ be tuples in $\mathcal{A}$.
\begin{enumerate}

\item  We say that $\overline{a}\equiv^0\overline{b}$ if $\overline{a}$ and $\overline{b}$ satisfy the same quantifier-free formulas,

\item  For $\alpha > 0$, we say that $\overline{a}\equiv^\alpha\overline{b}$ if for all $\beta < \alpha$, for each $\overline{c}$, there exists $\overline{d}$, and for each $\overline{d}$, there exists $\overline{c}$, such that $\overline{a},\overline{c}\equiv^\beta\overline{b},\overline{d}$.

\end{enumerate}

\end{defn}

\begin{defn}\ 

\begin{enumerate}

\item  The \emph{Scott rank of a tuple $\overline{a}$} in
  $\mathcal{A}$ is the least $\beta$ such that for all $\overline{b}$,
  the relation $\overline{a}\equiv^\beta\overline{b}$ implies
  $(\mathcal{A},\overline{a})\cong (\mathcal{A},\overline{b})$.

\item  The \emph{Scott rank of $\mathcal{A}$}, $SR(\mathcal{A})$, is the least ordinal $\alpha$ greater than the ranks of all tuples in $\mathcal{A}$.

\end{enumerate}

\end{defn}

\noindent
\textbf{Example}:  If $\mathcal{A}$ is an ordering of type $\omega$, then $SR(\mathcal{A}) = 2$.
We have $\overline{a}\equiv^0\overline{b}$ if $\overline{a}$ and $\overline{b}$ are ordered in the same way.  We have $\overline{a}\equiv^1\overline{b}$ if the corresponding intervals (before the first element and between successive elements) have the same size, and this is enough to assure isomorphism.  From this, it follows that the tuples have Scott rank $1$, so the ordering itself has Scott rank $2$.    

\bigskip

We are interested in \emph{computable} structures.  We adopt the following conventions.

\begin{enumerate}

\item  Languages are computable, and each structure has for its universe a subset of $\omega$. 

\item  We identify a structure $\mathcal{A}$ with its atomic diagram $D(\mathcal{A})$.  

\item  We identify sentences with their G\"{o}del numbers.

\end{enumerate}

By these conventions, a structure $\mathcal{A}$ is \emph{computable} (or \emph{arithmetical}) if $D(\mathcal{A})$, thought of as a subset of $\omega$, is computable (or arithmetical).  

\bigskip

\emph{Computable infinitary formulas} are useful in describing computable structures.  Roughly speaking, these are infinitary formulas in which the disjunctions and conjunctions are over c.e.\ sets.  They are essentially the same as the formulas in the least admissible fragment of $L_{\omega_1\omega}$.  For a more precise description of computable infinitary formulas, see 
\cite{A-K}.        

We may classify computable infinitary formulas as \emph{computable $\Sigma_\alpha$}, or \emph{computable $\Pi_\alpha$}, for various computable ordinals $\alpha$.  We have the useful 
fact that in a computable structure, a relation defined by a computable $\Sigma_\alpha$ (or 
computable $\Pi_\alpha$) formula will be $\Sigma^0_\alpha$ (or $\Pi^0_\alpha$).  
To illustrate the expressive power of computable infinitary formulas, we note that there is a 
natural computable $\Pi_2$ sentence characterizing the class of Abelian $p$-groups.  For 
each computable ordinal $\alpha$ there is a computable $\Pi_{2\alpha}$ formula saying (of 
an element of an Abelian $p$-group), that the height is at least $\omega\cdot\alpha$. 

We have a version of compactness for computable infinitary formulas.  

\begin{thm} [Barwise--Kreisel Compactness] 
\label{thm1.2} 

Let $\Gamma$ be a $\Pi^1_1$ set of computable infinitary sentences.  If every $\Delta^1_1$ subset of $\Gamma$ has a model, then $\Gamma$ has a model.

\end{thm} 

Barwise-Kreisel Compactness differs from ordinary Compactness in that it can be used to produce computable structures.

\begin{cor}
\label{cor1.3}

Let $\Gamma$ be a $\Pi^1_1$ set of computable infinitary sentences.  If every $\Delta^1_1$ subset has a computable model, then $\Gamma$ has a computable model. 

\end{cor}

The next two corollaries give evidence of the expressive power of computable infinitary formulas.  

\begin{cor}
\label{cor1.4} 

If $\mathcal{A}$, $\mathcal{B}$ are computable structures satisfying the same computable infinitary sentences, then $\mathcal{A}\cong\mathcal{B}$. 

\end{cor}

\begin{cor}
\label{cor1.5}  

Suppose $\overline{a}$, $\overline{b}$ are tuples satisfying the same computable infinitary formulas in a computable structure $\mathcal{A}$.  Then there is an automorphism of $\mathcal{A}$ taking $\overline{a}$ to $\overline{b}$.

\end{cor} 

Corollary \ref{cor1.5} yields a bound on the Scott ranks for
computable structures \cite{N1}.

\begin{prop}
\label{prop1.6} 

For a computable structure $\mathcal{A}$, we have $SR(\mathcal{A})\leq\omega_1^{CK}+1$. 

\end{prop}

The Barwise-Kreisel Compactness Theorem and the three corollaries are all well known, and may be found in \cite{A-K}.  One point in the proof of the Barwise-Kreisel Compactness Theorem is expanded in \cite{G-K}.  The following observation is given in \cite{K-Y}, among other places.    

\begin{prop}
\label{prop1.7}

For a computable structure $\mathcal{A}$, 

\begin{enumerate}

\item  $SR(\mathcal{A}) < \omega_1^{CK}$ if there is some computable ordinal $\beta$ such that
the orbits of all tuples are defined by computable $\Pi_\beta$ formulas.

\item  $SR(\mathcal{A}) = \omega_1^{CK}$ if the orbits of all tuples are defined by computable infinitary formulas, but there is no computable bound on the complexity of these formulas.

\item  $SR(\mathcal{A}) = \omega_1^{CK}+1$ if there is some tuple whose orbit is not defined by any computable infinitary formula.

\end{enumerate}

\end{prop}

Low Scott rank is associated with simple Scott sentences.
A \emph{Scott sentence} for $\mathcal{A}$ is a sentence whose countable models are just the isomorphic copies of $\mathcal{A}$ (as in the Scott Isomorphism Theorem).
Nadel \cite{N1}, \cite{N2} showed the following.  

\begin{thm} [Nadel] 
\label{thm1.8}

For a computable structure $\mathcal{A}$, 
$SR(\mathcal{A})$ is computable iff $\mathcal{A}$ has a computable 
infinitary Scott sentence.    

\end{thm}

We turn to examples of computable structures illustrating the different possible Scott ranks.  
There are familiar examples of computable structures of computable rank.

\begin{prop}
\label{prop1.9} 

For the following classes of structures, all computable members have computable Scott rank:

\begin{enumerate}

\item  well orderings,

\item  superatomic Boolean algebras,

\item  reduced Abelian $p$-groups.

\end{enumerate}   

\end{prop}
 
There are some well-known examples of computable structures of Scott
rank $\omega_1^{CK}+1$.  Harrison \cite{H} showed that there is a
computable ordering of type $\omega_1^{CK}(1+\eta)$.  This ordering,
the \emph{Harrison ordering}, gives rise to some other computable
structures with similar properties.  The \emph{Harrison Boolean
  algebra} is the interval algebra of the Harrison ordering.  The
\emph{Harrison Abelian $p$-group} has length $\omega_1^{CK}$, with all
infinite Ulm invariants, and with a divisible part of infinite
dimension.

\begin{prop}
\label{prop1.10} 

The Harrison ordering, Harrison Boolean algebra, and Harrison Abelian $p$-groups all have Scott rank $\omega_1^{CK}+1$.

\end{prop}

For the Harrison ordering, the rank is witnessed by any element $a$
outside the initial copy of $\omega_1^{CK}$.  Similarly, in the
Harrison Boolean algebra, the rank is witnessed by any non-superatomic
element, and in the Harrison Abelian\linebreak $p$-group, the rank is
witnessed by any divisible element.

\bigskip

The Harrison ordering has further interesting features.  First, the
computable infinitary sentences true in the Harrison ordering are all
true in orderings of type $\omega_1^{CK}$, so the conjunction of these
sentences is not a Scott sentence.  Second, although there are many
automorphisms, there is at least one computable copy in which there is no
non-trivial hyperarithmetical automorphism.

\bigskip

For Scott rank $\omega_1^{CK}$, it is not so easy to find computable examples.  There is an arithmetical example in \cite{M}.    

\begin{thm} [Makkai] 
\label{thm1.11}

There is an arithmetical structure $\mathcal{A}$ of rank $\omega_1^{CK}$.

\end{thm}

For Makkai's example, in contrast to the Harrison ordering, the set of computable infinitary sentences true in the structure is $\aleph_0$ categorical, so the conjunction of these sentences is a Scott sentence for the structure.  The structure shares with the Harrison ordering, as originally constructed, the feature that although there are many automorphisms, there is no non-trivial hyperarithmetical automorphism.   In \cite{K-Y}, Makkai's result is refined as follows.  

\begin{thm} 
\label{thm1.12} 

There is a computable structure of Scott rank $\omega_1^{CK}$.   

\end{thm}
 
In the remainder of the present section, we will review the results of
\cite{K-Y} and \cite{C-K-Y} establishing Theorem \ref{thm1.12}, first
for abstract structures, then for trees.  In Section \ref{secexs} we
will prove new results showing that there are computable undirected
graphs, linear orderings, and fields with Scott rank $\omega_1^{CK}$.
In Section \ref{secxcl} we will demonstrate that there are no
computable Abelian $p$-groups and no computable models of either the
computable infinitary theory of well-orderings or that of superatomic
Boolean algebras with Scott rank $\omega_1^{CK}$.  Finally, in Section
\ref{secapprox} we show that the examples constructed in this paper
are strongly computably approximable.

\subsection{Known Computable Structures of Scott Rank $\omega_1^{CK}$}

In \cite{K-Y}, there are two different proofs of Theorem
\ref{thm1.12}.  The first takes Makkai's example and, without
examining it, codes it into a computable structure in a way that
preserves the rank.
The second is a re-working of Makkai's construction, which incorporates a suggestion of 
Shelah (given at the end of Makkai's paper), and a suggestion of Sacks.  The structure is 
a ``group tree'' $\mathcal{A}(\mathcal{T})$, derived from a tree $\mathcal{T}$.  Morozov \cite{Mo} 
used the same construction.  He showed that if $\mathcal{T}$ is a
computable tree having a path but no hyperarithmetical path, then
$\mathcal{A}(\mathcal{T})$ is a computable structure which has the feature of having many
automorphisms but no non-trivial automorphism.  The Harrison ordering
shares this feature.
To get a $\mathcal{A}(\mathcal{T})$ as in Theorem \ref{thm1.12}, we need a tree $\mathcal{T}$ with special properties.  We need some definitions to state these properties.  Let $\mathcal{T}$ be a subtree of $\omega^{<\omega}$.  We define \emph{tree rank} for elements of $\mathcal{T}$, and for $\mathcal{T}$ itself.   

\begin{defn}\

\begin{enumerate}

\item  $rk(\sigma) = 0$ if $\sigma$ is terminal, 

\item  for $\alpha > 0$, $rk(\sigma) = \alpha$ if all successors of $\sigma$ have ordinal rank, and $\alpha$ is the first ordinal greater than these ordinals,

\item  $rk(\sigma) = \infty$ if $\sigma$ does not have ordinal rank.

\end{enumerate}
We let $rk(\mathcal{T}) = rk(\emptyset)$.

\end{defn}

\noindent
\textbf{Fact}.  $rk(\sigma) = \infty$ iff $\sigma$ extends to a path.  

\bigskip

For a tree $\mathcal{T}$, we let $\mathcal{T}_n$ be the set of elements at level $n$ in the tree---$\mathcal{T}_n = \mathcal{T}\cap\omega^n$.  

\begin{defn}

The tree $\mathcal{T}$ is \emph{thin} provided that for all $n$, the set of ordinal ranks of elements of $\mathcal{T}_n$ has order type at most $\omega\cdot n$.

\end{defn}

The following fact explains the importance of thinness.  

\bigskip
\noindent
\textbf{Fact}:  If $\mathcal{T}$ is a computable thin tree, then for each $n$, there is some computable $\alpha_n$ such that for all $\sigma\in\mathcal{T}_n$, if $rk(\sigma)\geq\alpha_n$, then 
$rk(\sigma) = \infty$. 

\bigskip

In \cite{K-Y}, we show the following.     
 
\begin{thm}\
\label{thm1.13}

\begin{enumerate}

\item  There exists a computable thin tree $\mathcal{T}$ with a path but no hyperarithmetical path.

\item  If $\mathcal{T}$ is a computable thin tree with a path but no hyperarithmetical path, then $\mathcal{A}(\mathcal{T})$ is a computable structure of Scott rank $\omega_1^{CK}$.

\end{enumerate}

\end{thm}

In \cite{C-K-Y}, we show that there is a computable tree of Scott rank $\omega_1^{CK}$.  
The idea is to take trees as in \cite{K-Y} and add a homogeneity property.  

\begin{defn} 

A tree $\mathcal{T}$ is \emph{rank-homogeneous} provided that for all $n$, 

\begin{enumerate}

\item  for all $\sigma\in\mathcal{T}_n$ and all computable $\alpha$, if there exists $\tau\in\mathcal{T}_{n+1}$ such that 
$rk(\tau) = \alpha < rk(\sigma)$, then $\sigma$ has infinitely many successors $\sigma'$ with\linebreak 
$rk(\sigma')=\alpha$. 

\item  for all $\sigma\in\mathcal{T}_n$, if $rk(\sigma) = \infty$, then $\sigma$ has infinitely many successors $\sigma'$ with $rk(\sigma')=\infty$. 

\end{enumerate}

\end{defn}

\noindent
\textbf{Fact}.  If $\mathcal{T}$ and $\mathcal{T}'$ are
rank-homogeneous trees, and for all $n$ there is an element in $T_n$
of rank $\alpha \in Ord\ \cup \{\infty\}$ if and only if there is an
element in $T'_n$ of rank $\alpha$, then
$\mathcal{T}\cong\mathcal{T}'$.

\bigskip 

In \cite{C-K-Y}, we obtain a tree of rank $\omega_1^{CK}$ as follows.  

\begin{thm} 
\label{thm1.15}\

\begin{enumerate}

\item  There is a computable, thin, rank-homogeneous tree $\mathcal{T}$ such that 
$rk(\mathcal{T})=~\infty$ but $\mathcal{T}$ has no hyperarithmetical path.

\item  If $\mathcal{T}$ is a computable, thin, rank-homogeneous tree such that 
$rk(\mathcal{T}) = \infty$ but $\mathcal{T}$ has no hyperarithmetical path, then $SR(\mathcal{T}) = \omega_1^{CK}$.

\end{enumerate}

\end{thm}

Like the group-trees, the trees in \cite{C-K-Y} have the feature that the computable infinitary theory is $\aleph_0$ categorical.  Unlike the group-trees, these trees have many non-trivial hyperarithmetical automorphisms.  

\begin{prop}
\label{prop1.15}

Suppose $\mathcal{T}$ is a computable rank-homogeneous tree.  Then

\begin{enumerate}

\item   $SR(\mathcal{T}) < \omega_1^{CK}$ if there is a computable bound on the ordinal tree ranks that occur in $\mathcal{T}$,

\item  $SR(\mathcal{T}) = \omega_1^{CK}$ if for each $n$, there is a computable bound on the ordinal tree ranks that occur in $\mathcal{T}_n$, but there is no computable bound on the ordinal tree ranks in $\mathcal{T}$,

\item  $SR(\mathcal{T}) = \omega_1^{CK}+1$ if there is some $n$ such that there is no computable bound on the ordinal tree ranks of tuples in $\mathcal{T}_n$.  

\end{enumerate}   

\end{prop}  

\section{Further examples}\label{secexs}

In this section, we give new examples of computable structures of Scott rank~$\omega_1^{CK}$.  

\begin{thm}
\label{thm2.1}  

Each of the following classes contains computable structures of Scott rank $\omega_1^{CK}$:

\begin{enumerate}

\item undirected graphs

\item  linear orderings,

\item  fields of any characteristic.

\end{enumerate} 

\end{thm}

\subsection{Computable embeddings}

We shall use a kind of computable embedding defined in \cite{C-C-K-M}.  
Let $K$ and $K'$ be classes of structures.  We suppose that each
structure has universe a subset of $\omega$.  Each class consists of
structures for a fixed computable language.  Moreover, each class is
closed under isomorphism, modulo the restriction on the universes.
Let $\Phi$ be a c.e.\ set $\Phi$ of pairs $(\alpha,\varphi)$, where
$\alpha$ is a finite set appropriate to be a subset of the atomic
diagram of a structure in $K$, and $\varphi$ is a sentence appropriate
to be in the atomic diagram of a structure in $K'$.  For each
$\mathcal{A}\in K$, let $\Phi(\mathcal{A})$ be the set of $\varphi$
such that for some $\alpha\subseteq D(\mathcal{A})$,
$(\alpha,\varphi)\in\Phi$.  Suppose for all $\mathcal{A}\in K$, the set
$\Phi(\mathcal{A})$ is the atomic diagram of some $\mathcal{B}\in K'$.
We identify the structure with its atomic diagram.  Now, $\Phi$ is a
\emph{computable embedding} of $K$ in $K'$ if for all
$\mathcal{A},\mathcal{A}'\in K$, we have $\mathcal{A}\cong\mathcal{A}'$ iff
$\Phi(\mathcal{A})\cong\Phi(\mathcal{A})$. \

\bigskip
\noindent
\textbf{Remark}:  If $\Phi$ is a computable embedding of $K$ in $K'$, and $\mathcal{A}$ is a computable structure in $K$, then $\Phi(\mathcal{A})$ is a member of $K'$ with a computable copy.

\bigskip

The following result is proved in \cite{K}.  

\begin{thm}
\label{thm2.2}

If $\Phi$ is a computable embedding of $K$ in $K'$, then for any computable infinitary formula $\varphi$, we can find a computable infinitary formula $\varphi^*$ such that for all $\mathcal{A}$ in $K$, $\Phi(\mathcal{A})\models\varphi$ iff $\mathcal{A}\models\varphi^*$.  Moreover, $\varphi^*$ has the same complexity as $\varphi$; i.e., if $\varphi$ is computable $\Sigma_\alpha$, then so is $\varphi^*$.    

\end{thm}

Using Theorem \ref{thm2.2}, we get the following.  

\begin{cor}
\label{cor2.3}

Let $\Phi$ be a computable embedding of $K$ in $K'$, where $K$ is axiomatized by a computable infinitary sentence.  For a hyperarithmetical structure $\mathcal{A}\in K$, if $SR(\Phi(\mathcal{A}))$ is computable, then $SR(\mathcal{A})$ is also computable.  

\end{cor}

\begin{proof}

Suppose $SR(\Phi(\mathcal{A}))$ is computable.  By Nadel's Theorem
(Theorem \ref{thm1.8}), $\Phi(\mathcal{A})$ has a computable
infinitary Scott sentence $\varphi$.  Let $\varphi^*$ be as guaranteed
by Theorem~\ref{thm2.2}.  If $\psi$ is a computable infinitary sentence
axiomatizing $K$, then $\psi\ \&\ \varphi^*$ is a Scott sentence for
$\mathcal{A}$.  Then, by Nadel's Theorem again, $SR(\mathcal{A})$ is
computable.

\end{proof}     

To prove Theorem \ref{thm2.1}, we shall describe a computable embedding of trees in undirected graphs, one of undirected graphs in fields of the desired characteristic, and one of undirected graphs in linear orderings.  We show that each of these embeddings has the following property.   

\bigskip
\noindent
\textbf{Definition}.  
Let $\Phi$ be a computable embedding of $K$ into $K'$.  We say that $\Phi$ has the
 \emph{rank-preservation property} provided that for all computable $\mathcal{A}\in K$, and $\mathcal{B} = \Phi(\mathcal{A})$, either $SR(\mathcal{A})$, $SR(\mathcal{B})$ are both computable, or else they are equal.
 
\bigskip

Corollary \ref{cor2.3} says that for a computable embedding $\Phi$
from $K$ into $K'$ and a computable structure $\mathcal{A}\in K$, if
$SR(\mathcal{A})\geq\omega_1^{CK}$, then
$SR(\Phi(\mathcal{A}))\geq\omega_1^{CK}$.  For rank preservation, we
need more.  In particular, we must show that if the orbits in
$\mathcal{A}$ are hyperarithmetical, then so are the orbits in
$\Phi(\mathcal{A})$, and if there is a bound on the complexity of the
orbits in $\mathcal{A}$ (i.e., all are $\Delta^0_\alpha$, for some
computable ordinal $\alpha$), then there is a bound on the complexity
of the orbits in $\Phi(\mathcal{A})$.  Actually, while this second
point follows from Corollary \ref{cor2.3}, we shall prove rank
preservation directly, without appealing to Corollary \ref{cor2.3}.
   
\bigskip

Supposing that we have the desired computable embeddings, with the
rank preservation property, we obtain Theorem \ref{thm2.1} as follows.
Let $\mathcal{T}$ be a computable tree of  Scott rank $\omega_1^{CK}$.
We get a computable undirected graph of Scott rank $\omega_1^{CK}$ by
first taking the image of $\mathcal{T}$ under under a rank-preserving
computable embedding of trees into undirected graphs and then passing
to a computable copy $\mathcal{G}$.  In the same way, we obtain
examples of fields and linear orderings.

\bigskip             

The following result of Soskov  \cite{So}, which is re-worked in \cite{G-H-K-S}, will be useful in calculating complexities of orbits.

\begin{thm} [Soskov]
\label{thm2.4}

Suppose $\mathcal{A}$ is a hyperarithmetical structure, and let $R$ be a relation on $\mathcal{A}$.  If $R$ is invariant under automorphisms, and hyperarithmetical, then it is definable in $\mathcal{A}$ by a computable infinitary formula.  

\end{thm}

This result implies that if an invariant relation $R$ is hyperarithmetical in one hyperarithmetical copy of a given structure $\mathcal{A}$, then in all hyperarithmetical copies of $\mathcal{A}$, the image of $R$ is hyperarithmetical.

\begin{cor}
\label{cor2.5}

Suppose $\mathcal{A}$ is hyperarithmetical.  Then $SR(\mathcal{A})$ is computable iff all of the orbits in $\mathcal{A}$ are hyperarithmetical, with a bound on the complexity (i.e., all are $\Delta^0_\alpha$, for some computable ordinal $\alpha$).    

\end{cor}

\begin{proof}

First, suppose $SR(\mathcal{A})$ is computable.  Then the orbits of tuples in $\mathcal{A}$ are defined by computable infinitary formulas of bounded complexity, so the orbits are all hyperarithmetical, with a bound on the complexity.  Now, suppose the orbits are all hyperarithmetical, with a bound on the complexity.  The \emph{orbit equivalence relation} is the relation that holds between a pair of tuples iff they are in the same orbit.  Let $\mathcal{A}^*$ be the variant of $\mathcal{A}$ with added elements representing the tuples from $\mathcal{A}$.  We include disjoint unary predicates $U_n$ representing $n$-tuples from $\mathcal{A}$, and we identify $U_1$ with the universe of $\mathcal{A}$, and put on this set the relations of $\mathcal{A}$.  For $n\geq 2$, we have projection functions $p^n_i$, for $1\leq i\leq n$, mapping each element of $U_n$ to the $i^{th}$ element of the corresponding tuple in $U_1$.  Clearly, $\mathcal{A}^*$ is hyperarithmetical, and the orbit equivalence relation is an invariant, hyperarithmetical relation on pairs in $\mathcal{A}^*$.  
By Theorem  \ref{thm2.4}, it is definable in $\mathcal{A}^*$ by a computable infinitary formula.  

\bigskip
\noindent
\textbf{Claim}:  There is a fixed $\alpha$ such that for all tuples $\overline{a}$ in $\mathcal{A}$, the orbit of $\overline{a}$ is defined by the conjunction of all computable $\Pi_\alpha$ formulas true of $\overline{a}$ in $\mathcal{A}$.  

\bigskip
\noindent
Proof of Claim:  Suppose not.  Let $\Gamma(x,y)$ be a $\Pi^1_1$ set of computable infinitary formulas saying that $x$ and $y$ are in different orbits, but they are in the same $U_n$ and for each $n$ and each computable infinitary formula $\varphi$ in variables $u_1,\ldots,u_n$, if $x,y\in U_n$, then $\varphi$ is satisfied by the tuple represented by $x$ iff it is satisfied by the tuple represented by $y$.  If there is no $\alpha$ as in the claim, then every hyperarithmetical subset of $\Gamma$ is satisfied by some pair in $\mathcal{A}^*$.  Therefore, the whole of $\Gamma$ is satisfied, a contradiction.

\bigskip

Using the claim, we get a bound on the complexity of formulas defining the orbits in $\mathcal{A}$.  Therefore, $\mathcal{A}$ has computable Scott rank.  

\end{proof}                      

\subsection{Embedding trees in undirected graphs}

There are several well-known methods for coding a tree in an undirected graph (see, for example, Marker \cite{Ma}).  We may represent a tree element $a$ by a point $r(a)$ with an edge connecting it to a triangle graph.  For a pair of tree elements $a,a'$, to indicate that $a'$ is a successor of $a$, we add a point $q(a,a')$, connected by an edge to a square, and we connect $r(a)$ and $r(a')$ to $q(a,a')$ by chains of length $2$, $3$, respectively.  All of these elements are distinct.  For convenience, we consider the top node of the tree to be a successor of itself.     

In \cite{C-C-K-M}, this idea is turned into a computable embedding.  We start with a large computable graph $\mathcal{G}$ including a representative $r(n)$, and attached triangle, for each $n\in\omega$, and also including a point $s(m,n)$, and attached square, allowing for the possibility that $n$ might be a successor of $m$.  For each tree $\mathcal{T}$, $\Phi(\mathcal{T})$ is the subgraph of $\mathcal{G}$ representing just the elements $n$ that are actually in $\mathcal{T}$ and the pairs $(m,n)$ that are actually in the successor relation in $\mathcal{T}$.  To show that the embedding has the rank preservation property, we note that there are finitary existential formulas $u(x)$ and $s(x,y)$ such that for any tree $\mathcal{T}$, $u$ and $s$ define in $\Phi(\mathcal{T})$ the universe and successor relation of a copy of $\mathcal{T}$.  For a computable tree $\mathcal{T}$, if $\mathcal{B} = \Phi(\mathcal{T})$ and $\mathcal{A}$ is the copy of $\mathcal{T}$ defined in $\mathcal{B}$ by the formulas $u$ and $s$, then we can see that $\mathcal{A}$ and $\mathcal{B}$ satisfy the hypotheses of the following proposition.    

\begin{prop}
\label{prop2.6}

Let $\mathcal{B}$ be a hyperarithmetical structure.  Suppose $\mathcal{A}$ is definable in $\mathcal{B}$ by computable infinitary formulas, and in case the language of $\mathcal{A}$ is infinite, there is a bound on the complexity of these formulas.  Suppose that all automorphisms of $\mathcal{A}$ extend to automorphisms of $\mathcal{B}$.  Finally, suppose that for each tuple $\overline{b}$ in $\mathcal{B}$, the orbit of $\overline{b}$ under automorphisms of $\mathcal{B}$ that fix $\mathcal{A}$ pointwise is definable by a computable infinitary formula $\psi(\overline{a},\overline{x})$, and there is a bound on the complexity of these formulas.  Then either $SR(\mathcal{A})$ and $SR(\mathcal{B})$ are both computable, or else they are equal.    

\end{prop}

\begin{proof}

Let $\overline{b}$ be a tuple in $\mathcal{B}$.  Let $\psi(\overline{a},\overline{x})$ 
define the orbit of $\overline{b}$ under automorphisms of $\mathcal{B}$ that fix the elements of $\mathcal{A}$.  Then $\overline{b}'$ is in the orbit of $\overline{b}$ in $\mathcal{B}$ iff there exists $\overline{a}'$ such that $\overline{a}'$ is in the orbit of $\overline{a}$ in $\mathcal{A}$ and $\mathcal{B}\models \psi(\overline{a}',\overline{b}')$.  Therefore, if the orbit of $\overline{a}$ in $\mathcal{A}$ is hyperarithmetical, so is the orbit of $\overline{b}$ in $\mathcal{B}$.  Moreover, if the orbits in $\mathcal{A}$ have bounded complexity, so do the orbits in $\mathcal{B}$.          
From this, it is clear that if $SR(\mathcal{A})$ is computable, so is $SR(\mathcal{B})$.    
If $SR(\mathcal{A}) = \omega_1^{CK}+1$, then there is some tuple $\overline{a}$ whose orbit is not defined by any computable infinitary formula.  By Soskov's Theorem, the orbit is not hyperarithmetical.  The orbit of $\overline{a}$ in $\mathcal{B}$ is the same, so $SR(\mathcal{B}) = \omega_1^{CK}+1$.  
Finally, suppose $SR(\mathcal{A}) = \omega_1^{CK}$.  The argument above shows that the orbits in $\mathcal{B}$ are all hyperarithmetical, since those in $\mathcal{A}$ are.  There is no bound on the complexity, since the orbits in $\mathcal{A}$ are among the orbits in $\mathcal{B}$.           

\end{proof}

\begin{cor}
\label{cor2.7}

There is a computable embedding $\Phi$ of trees into graphs such that $\Phi$ has the rank preservation property.  

\end{cor}

\subsection{Fields}

We obtain a computable embedding $\Phi$ of undirected graphs into
fields of any desired characteristic by modifying an embedding due to
Friedman and Stanley \cite{F-S}.  We describe the construction for
characteristic $\not= 2$.  Let  $\mathcal{F}$ be a computable
algebraically closed field with a computable sequence
$(b_n)_{n\in\omega}$ of algebraically independent elements, and such
that we can effectively determine the dependence relations.  For a
graph $\mathcal{G}$, the first step toward forming $\Phi(\mathcal{G})$
is to define the field $F_0$.  Let $F_-$ be a prime field of the
appropriate characteristic.  The field $F_0$ is the composite of all
the fields $acl(F_-(b_n))$.  We now form the field $\Phi(\mathcal{G})$
by adjoining the elements $\sqrt{c_i+c_j}$, where $i$ and $j$ are connected by an edge
in $\mathcal{G}$ and $c_i$ is inter-algebraic with $b_i$.  (For
characteristic $2$, the construction is similar except that we would
use cube roots instead of square roots.)

In the Friedman and Stanley embedding, the only added square roots were $\sqrt{b_i+b_j}$, where there is an edge connecting $i$ and $j$.  In \cite{C-C-K-M}, we observed that this gives a computable embedding.  The proof that the embedding preserves isomorphism is the same for the Friedman and Stanley embedding and the variant described above.  We need the fact that for all $d$ in $\Phi(\mathcal{G})$, if the algebraic closure of $d$ is present in $\Phi(\mathcal{G})$, then $d$ is interalgebraic with $b_i$ for some $i\in\mathcal{G}$.  We also need the fact that for $i,j\in\mathcal{G}$, not connected by an edge, there is no square root for $b_i+b_j$ in $\Phi(\mathcal{G})$. 

We must show that our computable embedding has the rank preservation property.  Let 
$\mathcal{G}$ be a computable graph.  If $\mathcal{B} = \Phi(\mathcal{G})$, and 
$\mathcal{A}$ is the copy of $\mathcal{G}$ with universe consisting of the algebraic closures of the special basis elements $b_i$, for $i\in\mathcal{G}$ and edge relation defined in terms of existence of square roots (or cube roots).  It is not difficult to see that $\mathcal{A}$ and $\mathcal{B}$ satisfy the conditions for the following.

\begin{prop}
\label{prop2.8}

Let $\mathcal{B}$ be a hyperarithmetical structure, and let
$\mathcal{A}$ be a definable quotient in $\mathcal{B}$; i.e., there
exist a structure $\mathcal{A}^* = (D,(R_i)_{i\in I})$ and a
congruence relation $\equiv$ such that $\mathcal{A}^*$ and $\equiv$
are definable in $\mathcal{B}$ by computable infinitary formulas of
bounded complexity, and $\mathcal{A} = \mathcal{A}^*/_\equiv$.
Further suppose that for any choice function $c: \mathcal{A} \to
\mathcal{A^*}$ where $c(a/_\equiv) \in a/_\equiv$ we have
$\mathcal{A}\cong_c\mathcal{A}_c$.  Suppose in addition that the
following conditions are satisfied.

\begin{enumerate}

\item  For any automorphism $f$ of $\mathcal{A}$ and any choice
  function $c$, the automorphism of $\mathcal{A}_c$, given by $c\circ
  f\circ c^{-1}$, extends to an automorphism of $\mathcal{B}$.

\item  For any tuple $\overline{b}$ in $\mathcal{B}$, the orbit of
  $\overline{b}$ under automorphisms of $\mathcal{B}$ that fix $D$
  pointwise is defined by a computable infinitary formula,
  $\varphi(\overline{d},\overline{x})$, of bounded complexity, where
  for any choice function $c$, the parameters $\overline{d}$ may be
  chosen to be in $\mathcal{A}_c$.

\end{enumerate} 
Then either $\mathcal{A}$ and $\mathcal{B}$ have the same Scott rank, or else both have computable Scott rank.  

\end{prop} 

\begin{proof}

There is a hyperarithmetical choice function $c$.  We have a hyperarithmetical copy $\mathcal{A}_c$ of $\mathcal{A}$ such that $\mathcal{A}\cong_c\mathcal{A}_c$.  Suppose $\overline{a}$, $\overline{a}'$ are tuples in $\mathcal{A}_c$.  If $\overline{a}$ and $\overline{a}'$ are in the same orbit in $\mathcal{A}_c$, then by 1,  they are in the same orbit in $\mathcal{B}$.  Conversely, if $\overline{a}$ and $\overline{a}'$ are in the same orbit in $\mathcal{B}$, the automorphism $f$ of $\mathcal{B}$ taking $\overline{a}$ to $\overline{a}'$ restricts to an automorphism of $\mathcal{A}^*$ taking the equivalence class of $a_i$ to that of $a_i'$.  We get an induced automorphism $f_c$ of $\mathcal{A}_c$ taking $\overline{a}$ to $\overline{a}'$.  It follows that if $SR(\mathcal{B})$ is computable, or $\leq\omega_1^{CK}$, then the same is true of $SR(\mathcal{A})$.      

Let $\overline{b}$ be a tuple in $\mathcal{B}$.  Take $\varphi(\overline{d},\overline{x})$ as in 2, defining the orbit of $\overline{b}$ over $D$, where the parameters $\overline{d}$ are in $\mathcal{A}_c$.  

\bigskip
\noindent
\textbf{Claim}:  $\overline{b}'$ is in the orbit of $\mathcal{B}$ iff there exists $\overline{d}'$ in the orbit of $\overline{d}$ in $\mathcal{A}_c$ such that $\varphi(\overline{d}',\overline{b}')$ holds in $\mathcal{B}$.  

\bigskip
\noindent
Proof of Claim:  First, suppose $\overline{b}'$ is in the orbit of
$\overline{b}$.  If $f$ is an automorphism of $\mathcal{B}$ taking
$\overline{b}$ to $\overline{b}'$, then, as above, $f$ restricts to an
automorphism of $\mathcal{A}^*$, and we get an automorphism $f_c$ of
$\mathcal{A}_c$, taking $c(d_i)$ to $c(f(d_i))$.  While
$f(\overline{d})$ may not be in $\mathcal{A}_c$, $\overline{d}' =
f_c(\overline{d})$ is in $\mathcal{A}_c$.  By 1, there is an
automorphism of $\mathcal{B}$ extending $f_c$, and we have
$\varphi(\overline{d}',\overline{b}')$.
Now, suppose $\varphi(\overline{d}',\overline{b}')$ holds in $\mathcal{B}$, where $\overline{d}'$ is in the orbit of $\overline{d}$ in $\mathcal{A}_c$.  By 1, an automorphism of $\mathcal{A}_c$ mapping $\overline{d}'$ to $\overline{d}$ extends to an automorphism $f$ of $\mathcal{B}$.  Then $f$ maps $\overline{b}'$ to a tuple $\overline{b}''$ satisfying $\varphi(\overline{d},\overline{x})$, and this $\overline{b}''$ is in the orbit of $\overline{b}$.  This completes the proof of the claim.

\bigskip 

Using the claim, we can see that if the orbit of $\overline{d}$ in $\mathcal{A}_c$ is hyperarithmetical, then the orbit of $\overline{b}$ in $\mathcal{B}$ is also hyperarithmetical.  Moreover, if the orbits of tuples in $\mathcal{A}_c$ have bounded complexity, then the orbits in $\mathcal{B}$ also have bounded complexity.  Therefore, if $SR(\mathcal{A})$ is computable, or $\leq\omega_1^{CK}$, then so is $SR(\mathcal{B})$.  Putting together what we have shown, we get the fact that either $SR(\mathcal{A})$ and $SR(\mathcal{B})$ are both computable or else they are equal.             

\end{proof} 

\begin{cor}
\label{cor2.9}

There is a computable embedding $\Phi$ of undirected graphs into fields of any desired characteristic, such that $\Phi$ has the rank preservation property.   

\end{cor}              

\subsection{Linear orderings}

We have a computable embedding of undirected graphs in linear orderings.  Friedman and Stanley gave a Borel embedding \cite{F-S}, which can be made computable.  We first form a large ordering $\mathcal{L}$, the result of putting the lexicographic ordering on $Q^{<\omega}$.  Let 
$(t_n)_{n\in\omega}$ be a list of the atomic types for tuples in graphs, such that those with 
$m$ variables appear before those with $m+1$ variables.   Let 
$(Q_a)_{a\in\omega}$ be a computable partition of 
$Q$ into  dense subsets.  The sets $Q_0$ and $Q_1$ have special roles.  Let 
$\mathcal{G}$ be a graph.  Then $\Phi(\mathcal{G})$ is the sub-ordering of $\mathcal{L}$ with elements $q_1r_1q_2r_2,\ldots q_n r_n k\in Q^{<\omega}$ such that for some finite sequence $a_1,\ldots,a_n$, say of atomic type $t_m$ in $\mathcal{G}$, we have $q_i\in Q_{a_i}$, for $i < n$, $r_i\in Q_0$, $r_n\in Q_1$, and $k < m$.   

The authors are grateful to Desmond Cummins for a detailed proof (in
work related to his senior thesis), that this $\Phi$ really is a
computable embedding.  Here is a brief sketch of the proof.
Suppose $\mathcal{G}\cong_f\mathcal{G}'$.  To show that
$\Phi(\mathcal{G})\cong\Phi(\mathcal{G}')$, it is enough to show that
a certain set $\mathcal{F}$ of finite partial $1-1$ functions has the
back-and-forth property. For $b = q_1,r_1,\ldots,q_n,r_n,k$, where
$q_i\in Q_{a_i}$, let $g(b) = (a_1,\ldots,a_n)$.  Let
$p\in\mathcal{F}$ if $p$ maps $(b_1,\ldots,b_n)$ in
$\Phi(\mathcal{G})$ to $(b_1',\ldots,b_n')$ in $\Phi(\mathcal{G}')$,
where $f$ maps $g(b_i)$ to $g(b_i')$, $b_i < b_j$ iff $b_i' < b_j'$,
$b_i$ and $b_i'$ have the same last term (so they have the same
position in their maximal discrete intervals), and if $b_i$ and $b_j$
have a common initial segment of length $2m-1$, then so do $b_i'$ and
$b_j'$.  We say how to extend $p$, adding $b_{n+1}$ to the domain.
Take the greatest $m$ such that for some $i\leq n$, $b_{n+1}$ agrees
with $b_i$ on an initial segment of length $2m-1$.  We extend $p$,
mapping $b_{n+1}$ to some $b_{n+1}'$ agreeing with $b_i'$ on the
initial segment of length $2m-1$, such that $b_i < b_{n+1}$, or $b_i >
b_{n+1}$ iff $b_i' < b_{n+1}'$, or $b_i' > b_{n+1}'$ respectively,
then $b_{n+1}$ has the same last term as $b_{n+1}$, and $f$ maps
$g(b_{n+1})$ to $g(b_{n+1}')$.

Now, suppose $\Phi(\mathcal{G})\cong_f\Phi(\mathcal{G}')$.  
To show that $\mathcal{G}\cong\mathcal{G}'$, it is enough to show that a 
certain set of finite partial 
$1-1$ functions has the back-and-forth property.  Let 
$p\in \mathcal{F}$ if $p$ maps $(a_1,\ldots,a_n)$ in 
$\mathcal{G}$ to $(a_1',\ldots,a_n')$ in 
$\mathcal{G}'$, and for some\linebreak
$b = (q_1,r_1,\ldots,q_n,r_n,0)$ in 
$\Phi(\mathcal{G})$, we have 
$f(b) = (q_1',r_1',\ldots,q_n',r_n',0)$, where 
$a_i\in Q_{a_i}$ and $a_i'\in Q_{a_i'}$.  We say how to extend $p$, adding $a_{n+1}\in\mathcal{G}$ to the domain (adding an element to the range is symmetric).  Take $d$ agreeing with $b$ down to $q_n$, with further terms $r_n^*,q_{n+1},r_{n+1},0$, where $q_{n+1}\in Q_{a_{n+1}}$.  Then $f(d)$ will agree with $f(b)$ down to $q_n'$, with further terms $r_n^**,q_{n+1}',r_{n+1}',0$.  Say $q_{n+1}'\in Q_{a_{n+1}'}$.  We extend $p$ mapping $a_{n+1}$ to $a_{n+1}'$.   

\bigskip

We must show that $\Phi$ has the rank preservation property.  

\bigskip
\noindent
\textbf{Claim 1}:  There is a computable mapping $f$ taking tuples in $\mathcal{G}$ to elements of $\Phi(\mathcal{G})$, such that $\overline{a}$ and $\overline{a}'$ are in the same orbit in 
$\mathcal{G}$ iff $f(\overline{a})$ and $f(\overline{a}')$ are in the same orbit in $\Phi(\mathcal{G})$.          

\bigskip
\noindent
Proof of Claim 1:  For each tuple $\overline{a} = (a_1,\ldots,a_n)$ in $\mathcal{G}$, we let $f(\overline{a})$ be  the element $q_1r_1q_2r_2,\ldots q_n r_n 0$, where $q_i$ is first in $Q_{a_i}$, for $i < n$, $r_i$ is first in $Q_0$, and $r_n$ is first in $Q_1$.  Then $\overline{a}$ and $\overline{a}'$ are in the same orbit in $\mathcal{G}$ iff their $f$-images are in the same orbit in $\Phi(\mathcal{G})$.  

\bigskip
\noindent
\textbf{Claim 2}:  There is a definable set $X\subseteq\Phi(\mathcal{G})$ with a computable mapping $g$ from $X$ to tuples in $\mathcal{G}$, such that for $b,b'\in X$, $b$ and $b'$ are in the same orbit in $\Phi(\mathcal{G})$ iff $g(b)$ and $g(b')$ are in the same orbit in $\mathcal{G}$.

\bigskip
\noindent
Proof of Claim 2:   We let $X$ consist of the sequences in $\Phi(\mathcal{G})$ ending in $0$.  These are the left limit points. Suppose $b\in X$, say $b = q_1r_1q_2r_2,\ldots q_n r_n 0$, where $q_i\in Q_{a_i}$, for $i < n$, $r_i\in Q_0$, and $r_n\in Q_1$.  We let $g(b) = (a_1,\ldots,a_n)$.  
For $b,b'\in X$, $b$ and $b'$ are in the same orbit in $\Phi(\mathcal{G})$ iff $g(b)$ and $g(b')$ are in the same orbit in $\mathcal{G}$. 

\bigskip
\noindent
\textbf{Claim 3}:  For each tuple $\overline{b}$ in $\Phi(\mathcal{G})$, there is a tuple $\overline{d}$ in $X$, and a computable infinitary formula $\varphi(\overline{u},\overline{x})$ such that $\Phi(\mathcal{G})\models \varphi(\overline{d},\overline{b})$, and $\overline{b}'$ is in the orbit of $\overline{b}$ iff there exists $\overline{d}'$ in $X$ such that each $d_i$ in $\overline{d}$ is in the same orbit as the corresponding $d_i'$ in $\overline{d}'$, and  
$\Phi(\mathcal{G})\models\varphi(\overline{d}',\overline{b}')$.

\bigskip
\noindent
Proof of Claim 3:  Let $\overline{b} = (b_1,\ldots,b_r)$ be a tuple in $\Phi(\mathcal{G})$.  For each $b_i$, we let 
$d_i$ be the first element of the maximal discrete set containing $b_i$.  
From the size of the maximal discrete set, we can recover the length
of the tuple $g(d_i)$.  If $d_i$ agrees with $d_j$ on the first $2m-1$
terms, so that $g(d_i)$ and $g(d_j)$ agree on the first $m$ terms,
then the interval between $d_i$ and $d_j$ consists of elements
representing extensions of the same tuple in $\mathcal{G}$ of length
$m$.  Then the pair $(b_i,b_j)$ satisfies a formula $a_m(x,y)$ saying
that for all $z$ in the interval between $x$ and $y$, the maximal
discrete set containing $z$ has size representing a tuple of length
at least $m$.  Conversely, if $(b_i,b_j)$ satisfies the formula
$a_m(x,y)$, where $b_i$ and $b_j$ lie on different maximal discrete
sets and each represents a tuple from $\mathcal{G}$ of length at least
$m$, then the tuples agree on the first $m$ terms.

Suppose $\mathcal{A}$ is a computable graph, and let $\mathcal{B} = \Phi(\mathcal{A})$.  Let $X$, $f$, and $g$ be as described above.  Then by the arguments above, $\mathcal{A}$ and $\mathcal{B}$ satisfy the hypotheses of the following result.   

\begin{thm}
\label{thm2.10}

Let $\mathcal{A}$ and $\mathcal{B}$ be hyperarithmetical structures.  

\begin{enumerate}

\item  Suppose there is a hyperarithmetical map $f$ from tuples in $\mathcal{A}$ to tuples in $\mathcal{B}$ such that $\overline{a}$ and $\overline{a}'$ are in the same orbit in $\mathcal{A}$ iff $f(\overline{a})$ and $f(\overline{b})$ are in the same orbit in $\mathcal{B}$.  Then if $SR(\mathcal{B})$ is computable, so is $SR(\mathcal{A})$, and if $SR(\mathcal{B})\leq\omega_1^{CK}$, then $SR(\mathcal{A})\leq\omega_1^{CK}$.  

\item  Suppose $g$ is a hyperarithmetical map from a set 
$X$ of tuples in $\mathcal{B}$, invariant under automorphism, to tuples in 
$\mathcal{A}$, such that for 
$\overline{d},\overline{d}'\in X$, $\overline{d}$ and $\overline{d}'$ are in the same orbit in 
$\mathcal{B}$ iff $g(\overline{d})$ and 
$g(\overline{d}')$ are in the same orbit in $\mathcal{A}$.  Suppose further that for 
each tuple $\overline{b}$ in 
$\mathcal{B}$, there is a finite collection of tuples 
$\overline{d}_1,\ldots,\overline{d}_n$ in $X$, and for some $\beta <
\alpha$ there is a computable $\Sigma_\beta$ formula $\varphi$ which
is true of
$\overline{d}_1,\ldots,\overline{d}_n,\overline{b}$, such that for all
$\overline{b}'$ in $\mathcal{B}$, $\overline{b}$ and $\overline{b}'$
are in the same orbit iff there exist
$\overline{d}_1',\ldots,\overline{d}_n'$ in $X$ such that
$\overline{d}_i$ and $\overline{d}_i'$ are in the same orbit, and
$\varphi$ is satisfied by
$\overline{d}_1',\ldots,\overline{d}_n,\overline{b}'$ in
$\mathcal{B}$.  Then if $SR(\mathcal{A})$ is computable, or
$\leq\omega_1^{CK}$, so is $SR(\mathcal{B})$.

\end{enumerate}    

\end{thm}

\begin{proof}

For 1, suppose $f$ is $\Delta^0_\alpha$.  If the orbits in $\mathcal{B}$ are all $\Delta^0_\alpha$, 
then so are the orbits in $\mathcal{A}$.  If the orbits in $\mathcal{B}$ are all hyperarithmetical, but not necessarily of bounded complexity, then the orbits in $\mathcal{A}$ are also all hyperarithmetical.  Therefore, if $SR(\mathcal{B})$ is computable, or $\leq\omega_1^{CK}$, then the same is true of $SR(\mathcal{A})$.    

For 2, suppose $g$ is $\Delta^0_\alpha$.  If the orbits in $\mathcal{A}$ are all $\Delta^0_\alpha$, or all hyperarithmetical, then the same is true of the orbits of tuples in $X$.  Take a tuple $\overline{b}$ in $\mathcal{B}$, and let $\overline{d}_1,\overline{d}_n$ and $\varphi$ be as in 2.  Then the orbit of $\overline{b}$ is $\Delta^0_\alpha$, or hyperarithmetical, depending on the complexity of the orbits of certain tuples in $X$.  Therefore, if $SR(\mathcal{A})$ is computable, or $\leq\omega_1^{CK}$, the same is true of $SR(\mathcal{B})$.        

\end{proof}

\begin{cor}
\label{cor2.11}

There is a computable embedding $\Phi$ of graphs into linear orderings such that $\Phi$ has the rank preservation property.    

\end{cor}

\noindent
\textbf{Remarks}:  Part 2 of Theorem \ref{thm2.10} implies Part 1, with the roles of $\mathcal{A}$ and $\mathcal{B}$ reversed.  If $(\mathcal{A},\mathcal{B})$ satisfies 1 and 2, or $(\mathcal{A},\mathcal{B})$ and $(\mathcal{B},\mathcal{A})$ both satisfy 2, then either $\mathcal{A}$ and $\mathcal{B}$ both have computable Scott rank or else the ranks are the same.  We can also show that this implies our earlier general results, Proposition \ref{prop2.6} and Proposition \ref{prop2.8}.   

\bigskip
\noindent
\textbf{Proposition \ref{prop2.6}}.  
Suppose $\mathcal{A}$ is definable in $\mathcal{B}$ and every automorphism of $\mathcal{A}$ extends to an automorphism of $\mathcal{A}$.  Suppose for all $\overline{b}$ in $\mathcal{B}$, there is a formula $\varphi(\overline{c},\overline{b})$, of bounded complexity, defining the orbit of $\mathcal{B}$ over $\mathcal{A}$.  Then either $\mathcal{A}$ and $\mathcal{B}$ both have computable Scott rank or else the Scott ranks are the same. 

\bigskip
\noindent
Proof of Proposition \ref{prop2.6} from Theorem \ref{thm2.10}:  Let $f$ be the identity function on tuples from $\mathcal{A}$.  If $\overline{a}$ and $\overline{a}'$ are in the same orbit in $\mathcal{A}$, then they are in the same orbit in $\mathcal{B}$.  The converse is obvious.  Applying 1 above, we conclude that if $SR(\mathcal{B})$ is computable, or $\leq\omega_1^{CK}$, then $SR(\mathcal{A})$ is computable, or $\leq\omega_1^{CK}$.    

Let $g$ also be the identity function on tuples from 
$\mathcal{A}$.  Let $\overline{b}$ be a tuple in $\overline{B}$, and let 
$\varphi(\overline{c},\overline{x})$ define the orbit of $\mathcal{B}$ over 
$\mathcal{A}$, as in the hypothesis.  Then $\overline{b}'$ is in the same orbit as $\overline{b}$ iff there exists $\overline{c}'\in \mathcal{A}$, in the orbit of $\overline{c}$, such that $\varphi(\overline{c}',\overline{b}')$ holds.  Applying 2 above, we conclude that if $SR(\mathcal{A})$ is computable, or 
$\leq\omega_1^{CK}$, then $SR(\mathcal{A})$ is computable, or $\leq\omega_1^{CK}$.

\bigskip
\noindent
\textbf{Proposition \ref{prop2.8}}.  Let $\mathcal{B}$ be a
hyperarithmetical structure, and let $\mathcal{A} =
\mathcal{A}^*/_\equiv$, where $\equiv$ is a congruence relation on
$\mathcal{A}^* = (D,(R_i)_{i\in I})$, and $\mathcal{A}^*$ and $\equiv$
are definable in $\mathcal{B}$ by computable infinitary formulas of
bounded complexity.  Further suppose that for any choice function $c$ choosing one element
from each $a\in\mathcal{A}$, we have
$\mathcal{A}\cong_c\mathcal{A}_c$.  Suppose in addition that the following conditions
are satisfied.

\begin{enumerate}

\item  For any automorphism $f$ of $\mathcal{A}$ and any choice
  function $c$, the automorphism given by $c\circ f\circ c^{-1}$ of
  $\mathcal{A}_c$ extends to an automorphism of $\mathcal{B}$.

\item  For any tuple $\overline{b}$ in $\mathcal{B}$, the orbit of $\overline{b}$ under automorphisms of $\mathcal{B}$ that fix $D$ pointwise is defined by a computable infinitary formula, $\varphi(\overline{d},\overline{x})$, of bounded complexity, where the parameters $\overline{d}$ may be chosen to be in $ran(c)$ for any choice function $c$. 

\end{enumerate} 
Then either $\mathcal{A}$ and $\mathcal{B}$ have the same Scott rank, or else both have computable Scott rank.  

\bigskip
\noindent
Proof of Proposition \ref{prop2.8} using Theorem \ref{thm2.10}:  Let $c$ be a hyperarithmetical choice function.  We obtain a hyperarithmetical copy of $\mathcal{A}$ with universe equal to $ran(c)$.  We identify this with $\mathcal{A}$.  Let $f(\overline{a}) = \overline{a}$, for $\overline{a}$ in $\mathcal{A}$, and for $\overline{d}$ in $D$, let $g(\overline{d}) = \overline{a}$, where $c(d_i) = a_i$.  For any $\overline{b}$ in $\mathcal{B}$, we have a tuple $\overline{a}$ in $\mathcal{A}$ and a formula $\varphi(\overline{a},\overline{x})$ defining the orbit of $\overline{b}$ over $\mathcal{A}$.  Then $\overline{b}'$ is in the orbit of $\overline{b}$ iff there exists $\overline{d}$ in the orbit of $\overline{a}$ such that $\varphi(\overline{d},\overline{b}')$ holds.  

Take $\overline{a},\overline{a}'$ in $\mathcal{A}$.  If $\overline{a}$ and $\overline{a}'$ are in the same orbit in $\mathcal{A}$, then by 1 above, they are in the same orbit in $\mathcal{B}$.  Conversely, if they are in the same orbit in $\mathcal{B}$, then because $\mathcal{A}$ is a definable quotient, they are in the same orbit in $\mathcal{A}$.  Therefore, if $SR(\mathcal{B})$ is computable, or $\leq\omega_1^{CK}$, $SR(\mathcal{A})$ is also.  

Next, take $\overline{d}$, representing different equivalence classes in $D$.  Say $g(\overline{d}) = \overline{a}$ and $g(\overline{d}') = \overline{a}'$.  If $\overline{d}$ and $\overline{d}'$ are in the same orbit in $\mathcal{B}$, then $\overline{a}$ and $\overline{a}'$ are in the same orbit in $\mathcal{A}$, since $\mathcal{A}$ is a definable quotient structure.  If $\overline{a}$ and $\overline{a}'$ are in the same orbit in $\mathcal{A}$, then by 1 and 2 above, $\overline{d}$ and $\overline{d}'$ are in the same orbit in $\mathcal{B}$.  We have $\overline{d}'$ in the same orbit as $\overline{d}$ iff $g(\overline{d})$ is in the same orbit as $g(\overline{d}')$.  If $SR(\mathcal{A})$ is computable, or $\omega_1^{CK}$, then the orbits of tuples from $D$ have bounded complexity, or are all hyperarithmetical.   

Now, take $\overline{b}$ in $\mathcal{B}$, and let $\varphi(\overline{a},\overline{x})$ be as in 2 above, defining the orbit of $\overline{b}$ under automorphisms of $\mathcal{B}$ that fix $D$.  Then $\overline{b}'$ is in the orbit of $\overline{b}$ iff there exists $\overline{d}$ in the orbit of $\overline{a}$ such that $\varphi(\overline{d},\overline{b}')$ holds.  The complexity of the orbit of $\overline{b}$ is not far from the complexity of that of $\overline{a}$ (as a tuple in $D$).  Therefore, if $SR(\mathcal{A})$ is computable, or $\omega_1^{CK}$, then $SR(\mathcal{B})$ is also.                                                                        
     
\section{Classes with no computable member of Scott rank $\omega_1^{CK}$}\label{secxcl}

We have shown that there are computable structures of Scott rank
$\omega_1^{CK}$ in several familiar classes.  There are classes in
which there is not computable bound on the Scott ranks but there is no
a computable member of Scott rank $\omega_1^{CK}$.

\begin{prop}[essentially, Barwise]\label{noapg}

If $\mathcal{A}$ is a computable Abelian $p$-group, then $SR(\mathcal{A}) \not=\omega_1^{CK}$.   

\end{prop}

\begin{proof}

If $\mathcal{A}$ has computable length, then the Scott rank is
computable.  The only non-computable length possible for a computable
group is $\omega_1^{CK}$.  If $\mathcal{A}$ has length
$\omega_1^{CK}$, then it cannot be reduced.  The divisible elements
have Scott rank $\omega_1^{CK}$, so $\mathcal{A}$ has Scott rank
$\omega_1^{CK}+1$.

\end{proof}

The proofs of the following two results are essentially identical to
the proof of Proposition \ref{noapg}.

\begin{prop} If $\mathcal{A}$ is a model of the computable infinitary
  theory of well orderings, then $SR(\mathcal{A}) \neq \omega_1^{CK}$.
\end{prop}

\begin{prop} If $\mathcal{A}$ is a model of the computable infinitary
  theory of superatomic Boolean algebras, then $SR(\mathcal{A}) \neq
  \omega_1^{CK}$.
\end{prop}

\section{Strong computable approximability}\label{secapprox}

\begin{defn}

A structure $\mathcal{A}$ is \emph{strongly computably approximable} if for any $\Sigma^1_1$ set $S$, there is a uniformly computable sequence $(\mathcal{C}_n)_{n\in\omega}$ such that $n\in S$ iff $\mathcal{C}_n\cong\mathcal{A}$.  The structures $\mathcal{C}_n$ for $n\notin S$ are called \emph{approximating structures}.   

\end{defn}

For example, the Harrison ordering is strongly computably approximable by computable well orderings, where these all have computable Scott rank.  The following result is in \cite{C-K-Y}.

\begin{thm}
\label{thm4.1}

There is a computable tree $\mathcal{T}$, of Scott rank $\omega_1^{CK}$, such that $\mathcal{T}$ is strongly computably approximable.  Moreover, the approximating structures are trees of computable Scott rank.  

\end{thm}

The following result yields further structures that are strongly computably approximable.     

\begin{thm}
\label{thm4.2}

Suppose $\Phi$ is a computable embedding of $K$ into $K'$.  If $\mathcal{A}\in K$ is strongly computably approximable by structures in $K$, then $\Phi(\mathcal{A})$ is strongly computably approximable by structures in $K'$.  Moreover, if $\Phi$ has the rank preservation property, and the approximating structures for $\mathcal{A}$ have computable rank, then so do those for $\Phi(\mathcal{A})$.   

\end{thm}

\begin{proof}

If $\mathcal{A}$ is strongly computably approximable, then $SR(\mathcal{A})$ is not computable.  Therefore, $SR(\Phi(\mathcal{A}))$ is not computable.  Let $S$ be a $\Sigma^1_1$ set.  Take a uniformly computable sequence $(\mathcal{C}_n)_{n\in\omega}$ such that $\mathcal{C}_n\cong\mathcal{A}$ iff $n\in S$.  We get a uniformly computable sequence $(\mathcal{B}_n)_{n\in\omega}$, where $\mathcal{B}_n\cong\Phi(\mathcal{C}_n)$.  Then $\mathcal{B}_n\cong\Phi(\mathcal{A})$ iff $n\in S$.  Therefore, $\Phi(\mathcal{A})$ is strongly computably approximable by structures in $K'$.  Moreover, if $\Phi$ has the rank preservation property, and $\mathcal{C}_n$ has computable rank, so does $\Phi(\mathcal{C}_n)$.   

\end{proof}

Combining Theorem \ref{thm4.2} with the results in the previous section, we obtain the following.

\begin{thm}
\label{thm4.3}

In each of the following classes, there is a structure of Scott rank $\omega_1^{CK}$ that is strongly computably approximable by structures of computable Scott rank.

\begin{enumerate}

\item  undirected graphs

\item  linear orderings

\item  fields of any fixed characteristic

\end{enumerate}

\end{thm}

\begin{proof}

We have described a computable embedding of trees into undirected graphs, and computable embeddings of undirected graphs into linear orderings and fields of any desired characteristic.  These embeddings all have the rank preservation property.  Starting with a computable tree $\mathcal{T}$ of Scott rank $\omega_1^{CK}$, such that $\mathcal{T}$ is strongly computably approximable by trees of computable Scott rank, we obtain in each of the classes above, a structure $\mathcal{A}$ of Scott rank $\omega_1^{CK}$ such that $\mathcal{A}$ is strongly computably approximable by structures of computable Scott rank, in the given class.  

\end{proof} 

\section{Conclusion}

In the present paper, we have used computable embeddings to transfer results on trees to further classes of structures:  undirected graphs, fields of any desired characteristic, and linear orderings.  We used some known computable embeddings, modifying one of them slightly, and we showed that each has the rank preservation property.  Our results are not sensitive to the precise definition of computable embedding.  What we need is a function $\Phi$ from $K$ to $K'$ such that 

\begin{enumerate}

\item  for $\mathcal{A},\mathcal{A}'$ in $K$, $\mathcal{A}\cong\mathcal{A}'$ iff $\Phi(\mathcal{A})\cong\Phi(\mathcal{A}')$,

\item  if $\mathcal{A}\in K$ is computable, then $\Phi(\mathcal{A})$ has a computable copy $\mathcal{B}$, with index computable from that for $\mathcal{A}$,

\item  if $\mathcal{A}\in K$ is computable, then either $SR(\mathcal{A})$ and $SR(\Phi(\mathcal{A}))$ are both computable, or $SR(\mathcal{A}) = SR(\Phi(\mathcal{A}))$.

\end{enumerate}

Suppose $\Phi$ satisfies these three properties.  If $K$ contains a computable structure $\mathcal{A}$ of Scott rank $\omega_1^{CK}$, then $\Phi(\mathcal{A})$ is a structure in $K'$ of rank $\omega_1^{CK}$, with a computable copy.  Moreover, if $\mathcal{A}$ is strongly computably approximable, by structures in $K$ of computable Scott rank, then $\Phi(\mathcal{A})$ is strongly computably approximable, by structures in $K'$ of computable Scott rank.

\end{document}